\documentclass[11pt]{amsart}
\usepackage{amssymb}

\newtheorem{theorem}{Theorem}[section]
\newtheorem{corollary}[theorem]{Corollary}
\newtheorem{lemma}[theorem]{Lemma}
\newtheorem{proposition}[theorem]{Proposition}
\theoremstyle{definition}
\newtheorem{definition}[theorem]{Definition}

\newtheorem{question}[theorem]{Question}

\theoremstyle{remark}

\numberwithin{equation}{section}

\def\N{{\mathbb{N}}}

\def\3{{|\!|\!|}}

\begin{document}

\title[A hierarchy on non-archimedean CLI Polish groups]{A hierarchy on non-archimedean Polish groups admitting a compatible complete left-invariant metric}
\author{Longyun Ding}
\address{School of Mathematical Sciences and LPMC, Nankai University, Tianjin, 300071, P.R.China}
\email{dingly@nankai.edu.cn}
\thanks{Research is partially supported by the National Natural Science Foundation of China (Grant No. 11725103).}
\author{Xu Wang}
\email{1120210025@mail.nankai.edu.cn}



\begin{abstract}
    In this article, we introduce a hierarchy on the class of non-archimedean Polish groups that admit a compatible complete left-invariant metric. We denote this hierarchy by $\alpha$-CLI and L-$\alpha$-CLI where $\alpha$ is a countable ordinal. We establish three results:
    \begin{enumerate}
        \item $G$ is $0$-CLI iff $G=\{1_G\}$;
        \item $G$ is $1$-CLI iff $G$ admits a compatible complete two-sided invariant metric; and
        \item $G$ is L-$\alpha$-CLI iff $G$ is locally $\alpha$-CLI, i.e., $G$ contains an open subgroup that is $\alpha$-CLI.
    \end{enumerate}

    Subsequently, we show this hierarchy is proper by constructing non-archimedean CLI Polish groups $G_\alpha$ and $H_\alpha$ for $\alpha<\omega_1$, such that 
    \begin{enumerate}
        \item $H_\alpha$ is $\alpha$-CLI but not L-$\beta$-CLI for $\beta<\alpha$; and
        \item $G_\alpha$ is $(\alpha+1)$-CLI but not L-$\alpha$-CLI.
    \end{enumerate}
\end{abstract}
\maketitle

\section{Introduction}

A Polish group is {\it non-archimedean} if it has a neighborhood basis of its identity element consisting of open subgroups. By a theorem of Becker and Kechris (cf.~\cite[Theorem 1.5.1]{BK}), a Polish group is non-archimedean iff it is homeomorphic to a closed subgroup of $S_\infty$, the group of all permutations of $\N$ equipped with the pointwise convergence topology. A metric $d$ on a group $G$ is {\it left-invariant} if $d(gh,gk)=d(h,k)$ for all $g,h,k\in G$. A Polish group is {\it CLI} if it admits a compatible complete left-invariant metric.

Malicki \cite{malicki11} defined a notion of orbit tree $T_G$ for each closed subgroup $G$ of $S_\infty$, and showed that $G$ is CLI iff $T_G$ is well-founded. Moreover, he proved that the heights of orbit trees of all CLI closed subgroups of $S_\infty$ are cofinal in $\omega_1$. Malicki proved that the family of all CLI groups is coanalytic non-Borel based on this cofinality. After that, Xuan defined a different kind of orbit trees and showed that, a closed subgroup of $S_\infty$ is locally compact iff its orbit tree has finite height (cf. \cite[Theorem 3.7]{xuan}). It is worth noting that both kinds of orbit trees defined by Malicki and Xuan are all defined on closed subgroups of $S_\infty$ rather than directly on non-archimedean Polish groups. As a result, two topologically isomorphic closed subgroups of $S_\infty$ can have completely different orbit trees, and even the rank of their orbit trees can be different. This suggests that one cannot use ranks of orbit trees directly to define a hierarchy.

In this article, for a given non-archimedean CLI Polish group $G$, we use a neighborhood basis of the identity $1_G$ to define a new type of orbit trees. This differs from the approach which employs a closed subgroup of $S_\infty$ topologically isomorphic to $G$. More specifically, let $\mathcal G=(G_n)$ be a decreasing sequence of open subgroups of $G$ with $G_0=G$, such that $(G_n)$ forms a neighborhood basis of $1_G$. We will define a well-founded tree $T^{X(\mathcal G)}_\mathcal G$ and denote its rank by $\rho(\mathcal G)$. 
We let $\mathrm{rank}(G)$ be the ordinal given by $\max\{\beta:\rho(\mathcal{G})\geq\omega\cdot\beta\}$. Then we shall prove that the following are independent from the choice of $\mathcal{G}$: (a) the value of the ordinal $\textrm{rank}(G)$; and (b) whether $\rho(\mathcal{G})$ is a limit ordinal or not. These facts allow us to form a well-defined hierarchy on the class of non-archimedean CLI Polish groups: given an ordinal $\alpha<\omega_1$,
\begin{enumerate}
\item[(1)] if $\rho(\mathcal G)\le\omega\cdot\alpha$, we say $G$ is $\alpha$-CLI;
\item[(2)] if $\rho(\mathcal{G})\leq\omega\cdot\alpha+m$ for some $m<\omega$, i.e., ${\rm rank}(G)\le\alpha$, we say $G$ is L-$\alpha$-CLI.
\end{enumerate}
It is clear that, if $G$ is L-$\alpha$-CLI, then it is also $(\alpha+1)$-CLI.

The following theorem shows that the hierarchy classifies non-archimedean CLI Polish groups in a good manner:

\begin{theorem}
Let $G$ be a non-archimedean CLI Polish group and $\alpha$ be a countable ordinal. Then:
\begin{enumerate}
\item[(1)] $G$ is $0$-CLI iff $G=\{1_G\}$;
\item[(2)] $G$ is L-$0$-CLI iff $G$ is discrete;
\item[(3)] $G$ is $1$-CLI iff $G$ is TSI, i.e., $G$ admits a compatible complete two-sided invariant metric;
\item[(4)] $G$ is L-$\alpha$-CLI iff $G$ is locally $\alpha$-CLI, i.e., $G$ has an open subgroup which is $\alpha$-CLI.
\end{enumerate}
\end{theorem}

It is well-known that all compact Polish groups are TSI (c.f.~\cite[Theorem 2.1.5]{gaobook}), and all locally compact Polish groups are CLI (c.f.~\cite[Theorem 2.2.5]{gaobook}). Now we know that all compact non-archimedean Polish groups are $1$-CLI, and all locally compact non-archimedean Polish groups are L-$1$-CLI.

\begin{theorem}
Let $G$ be a non-archimedean CLI Polish group and $\alpha$ be a countable ordinal. Assume $H$ and $N$ are closed subgroups of $G$, and that $N$ is normal in $G$. If $G$ is $\alpha$-CLI (or L-$\alpha$-CLI), so are $H$ and $G/N$. In particular, we have ${\rm rank}(H)\le{\rm rank}(G)$ and ${\rm rank}(G/N)\le{\rm rank}(G)$.
\end{theorem}

\begin{theorem}
Let $(G^i)$ be a sequence of non-archimedean CLI Polish groups, $\alpha<\omega_1$, and let $G=\prod_iG^i$. Then we have
\begin{enumerate}
\item[(1)] $G$ is $\alpha$-CLI iff all $G^i$ are $\alpha$-CLI; and
\item[(2)] $G$ is L-$\alpha$-CLI iff all $G^i$ are L-$\alpha$-CLI and for all but finitely many $i$, $G^i$ is $\alpha$-CLI.
\end{enumerate}
\end{theorem}

Finally, we prove the following theorem, which indicates that this hierarchy is proper:

\begin{theorem}
For any $\alpha<\omega_1$, there exist two non-archimedean CLI Polish groups $G_\alpha$ and $H_\alpha$ with ${\rm rank}(G_\alpha)={\rm rank}(H_\alpha)=\alpha$ such that $H_\alpha$ is $\alpha$-CLI and $G_\alpha$ is L-$\alpha$-CLI but not $\alpha$-CLI.
\end{theorem}

\section{Preliminaries}

We denote the class of all ordinals by ${\rm Ord}$. For any $\alpha\in{\rm Ord}$, we define
$$\omega(\alpha)=\max\{0,\lambda:\lambda\le\alpha\mbox{ is a limit ordinal}\}.$$
Then $\alpha=\omega(\alpha)+m$ for some $m<\omega$.

Let $E$ be an equivalence relation on a set $X$, $x\in X$, and $A\subseteq X$. The $E$-equivalence class of $x$ is $[x]_E=\{y\in X:xEy\}$. Similarly, the $E$-saturation of $A$ is $[A]_E=\{y\in X:\exists z\in A\,(yEz)\}$.

The identity element of a group $G$ is denoted by $1_G$. Let $H$ be a subgroup of $G$, we denote the set of all left-cosets of $H$ by $G/H$.

A topological space is {\it Polish} if it is separable and completely metrizable. A topological group is {\it Polish} if its underlying topology is Polish. Let $G$ be a Polish group and $X$ a Polish space, an action of $G$ on $X$, denoted by $G\curvearrowright X$, is a map $a:G\times X\to X$ that satisfies $a(1_G,x)=x$ and $a(gh,x)=a(g,a(h,x))$ for $g,h\in G$ and $x\in X$. The pair $(X,a)$ is called a {\it Polish $G$-space} if $a$ is continuous. For brevity, we write $g\cdot x$ in place of $a(g,x)$. The {\it orbit equivalence relation} $E_G^X$ is defined as $xE_G^Xy\iff\exists g\in G\,(g\cdot x=y)$. Note that the $E_G^X$-equivalence class of $x$ is $G\cdot x=\{g\cdot x:g\in G\}$, which is also called the {\it $G$-orbit} of $x$. Similarly, for $A\subseteq X$, the $E_G^X$-saturation of $A$ is $G\cdot A=\{g\cdot x:g\in G\wedge x\in A\}$.

Let $<$ be a binary relation on a set $T$. We say that $(T,<)$ is a {\it tree} if
\begin{enumerate}
\item[(1)] $\forall s\in T\,(s\not<s)$,
\item[(2)] $\forall s,t,u\in T\,((s<t\wedge t<u)\Rightarrow s<u)$,
\item[(3)] $\forall s\in T\,(|\{t\in T:t<s\}|<\omega\wedge\forall t,u<s\,(t=u\vee t<u\vee u<t))$.
\end{enumerate}
For $s\in T$, we define ${\rm lh}(s)=|\{t\in T:t<s\}|$, which is called the length of $s$. It is clear that $s<t$ implies ${\rm lh}(s)<{\rm lh}(t)$. For $n<\omega$, we denote the $n$-th level of $T$ by
$$L_n(T)=\{s\in T:{\rm lh}(s)=n\}.$$
Each element in $L_0(T)$ is called a {\it root} of $T$.

Let $(T,<)$ be a tree. We say that $T$ is {\it well-founded} if any non-empty subset of $T$ contains at least a maximal element, or equivalently (under AC), if $T$ contains no infinite strictly increasing sequence. Let $T$ be a well-founded tree. We define the rank function $\rho_T:T\to{\rm Ord}$ by transfinite induction as
$$\rho_T(s)=\sup\{\rho_T(t)+1:s<t\wedge t\in T\}.$$
If $\rho_T(s)=0$, we say that $s$ is a {\it terminal} of $T$. Then we define
$$\rho(T)=\sup\{\rho_T(s)+1:s\in T\}.$$
So $\rho(T)=0$ iff $T=\emptyset$. It is clear that $\rho(T)=\sup\{\rho_T(s)+1:s\in L_0(T)\}$. If $L_0(T)=\{s_0\}$ is a singleton, then $\rho(T)=\rho_T(s_0)+1$ is a successor ordinal.

For $s\in T$, we define
$$T_s=\{t\in T:s=t\vee s<t\}.$$
Since $L_0(T_s)=\{s\}$, we have $\rho(T_s)=\rho_T(s)+1$. For the convenience of discussion, we let $T_s=\emptyset$ whenever $s\notin T$; in other words, $\rho(T_s)=0$. This convention will be useful in some proofs (see Lemma~\ref*{XtimesY}). Note that $\rho(T_s)$ is always a non-limit ordinal regardless of $s\in T$ or not.

First, we note the following facts:

\begin{proposition}\label{tree}
Let $T$ be a well-founded tree, then
$$\rho(T)=\sup\{\rho(T_s):s\in T\}=\sup\{\rho(T_s):s\in T\wedge s\in L_0(T)\},$$
and for all $s\in T$,
$$\begin{array}{ll}\rho(T_s) &=\sup\{\rho(T_t):s<t\wedge t\in T\}+1\cr
&=\sup\{\rho(T_t):s<t\wedge t\in T\wedge{\rm lh}(t)={\rm lh}(s)+1\}+1.\end{array}$$
\end{proposition}

\begin{proposition}\label{L_k(T)}
Let $(T,<)$ be a well-founded tree, $k<\omega$. Then we have
\begin{enumerate}
\item[(1)] $\sup\{\rho(T_s):s\in L_k(T)\}\le\rho(T)\le\sup\{\rho(T_s):s\in L_k(T)\}+k$;
\item[(2)] $\omega(\rho(T))=\omega(\sup\{\rho(T_s):s\in L_k(T)\})$;
\item[(3)] if $\rho(T_s)\ge\alpha$ for some $s\in L_k(T)$, then $\rho(T)\ge\alpha+k$.
\end{enumerate}
\end{proposition}

\begin{proof}
It is routine to prove clause (1) by induction on $k$ based on Proposition~\ref*{tree}. Clause (2) is an easy corollary of (1). And clause (3) is trivial.
\end{proof}

Let $(S,<)$ and $(T,<)$ be two trees. A map $\phi:S\to T$ is said to be an {\it order-preserving map} if
$$\forall s,t\in S\,(s<t\Rightarrow\phi(s)<\phi(t)).$$
It is said to be an {\it order-preserving embedding (isomorphism)} if it is injective (bijective) and
$$\forall s,t\in S\,(s<t\iff\phi(s)<\phi(t)).$$
In particular, an order-preserving map $\phi$ is said to be {\it Lipschitz} if ${\rm lh}(\phi(s))={\rm lh}(s)$ for all $s\in S$.

\begin{proposition}\label{subtree}
Let $(S,<)$ and $(T,<)$ be trees, and assume that $(T,<)$ is well-founded. If there exists an order-preserving map $\phi:S\to T$, then $(S,<)$ is also well-founded, and $\rho(S)\le\rho(T)$ holds.
\end{proposition}

\begin{proof}
If $S$ contains an infinite strictly increasing sequence $(s_n)$, then $(\phi(s_n))$ is an infinite strictly increasing sequence in $T$, contradicting that $(T,<)$ is well-founded.

Now we prove that $\rho_S(s)\le\rho_T(\phi(s))$ for $s\in S$ by induction on $\rho_S(s)$.
For the basis of induction, note that $\rho_S(s)=0\leq\rho_T(\phi(s))$ clearly holds; and for the inductive step, we have the following inequality:
$$\begin{array}{ll}\rho_S(s) &=\sup\{\rho_S(u)+1:s<u\wedge u\in S\}\cr
&\le\sup\{\rho_T(\phi(u))+1:s<u\wedge u\in S\}\cr
&\le\sup\{\rho_T(t)+1:\phi(s)<t\wedge t\in T\}=\rho_T(\phi(s)).\end{array}$$
This implies $\rho(S)\le\rho(T)$.
\end{proof}

\section{Definition of the hierarchy}

\begin{definition}
Let $X$ be a set, and $\mathcal E=(E_n)$ be a decreasing sequence of equivalence relations on $X$, i.e., $E_n\supseteq E_{n+1}$ for each $n<\omega$. We define
$$T_\mathcal E^X=\{(n,C):\exists x\in X\,(C=[x]_{E_n}\ne\{x\})\}.$$
For $(n,C),(m,D)\in T_\mathcal E^X$, we define
$$(n,C)<(m,D)\iff n<m\wedge C\supseteq D.$$
\end{definition}


It is straightforward to check that $(T_\mathcal E^X,<)$ is a tree. Note that $(n,C)\in T_\mathcal E^X$ iff $C$ is a non-singleton equivalence class of $E_n$. Here we omit all singleton classes in our definition. This will be crucial in the proof of Lemma~\ref*{CLItoWF}.

\begin{definition}
Let $G$ be a non-archimedean Polish group. We denote by ${\rm dgnb}(G)$ the set of all decreasing sequences $\mathcal G=(G_n)$ of open subgroups of $G$, such that $G_0=G$ and $(G_n)$ forms a neighborhood basis of $1_G$. Here `dgnb' stands for `decreasing group neighborhood basis'.

Let $X$ be a countable discrete Polish $G$-space, $\mathcal G=(G_n)\in{\rm dgnb}(G)$. We define $E_n=E_{G_n}^X$, i.e., $xE_ny\iff\exists g\in G_n\,(g\cdot x=y)$, and hence $[x]_{E_n}=G_n\cdot x$. Then we write $\mathcal E=(E_n)$ and
$$T_\mathcal G^X=T_\mathcal E^X.$$
\end{definition}

Therefore, $(n,C)\in T_\mathcal G^X$ iff $C$ is a non-singleton $G_n$-orbit. 
We point out that our definition is different from Malicki's and Xuan's, which are based on infinite orbits.

The following lemma builds a connection between non-archimedean CLI Polish groups and well-founded trees. Similar results also appear in \cite[Theorem 3]{malicki11} and \cite[Theorem 3.9]{xuan}.

\begin{lemma}\label{CLItoWF}
Let $G$ be a non-archimedean CLI Polish group, $X$ a countable discrete Polish $G$-space, and let $\mathcal G=(G_n)\in{\rm dgnb}(G)$. Then $T_\mathcal G^X$ is well-founded.
\end{lemma}

\begin{proof}
Assume for contradiction that $T_\mathcal G^X$ is ill-founded, then there exists an infinite sequence $(n,C_n),\,n<\omega$ in $T_\mathcal G^X$ with $C_n\supseteq C_{n+1}$ for each $n<\omega$.

Let $d$ be a compatible complete left-invariant metric on $G$.

Fix an $x_0\in C_0$. Then we have $G_0\cdot x_0=C_0\supseteq C_1$, so we can find a $g_0\in G_0$ such that $g_0\cdot x_0\in C_1$. Inductively, we can find a $g_n\in G_n$ for each $n<\omega$ such that $g_ng_{n-1}\cdots g_0\cdot x_0\in C_{n+1}$. Put $h_n=g_0^{-1}\cdots g_n^{-1}$ for each $n<\omega$. Then, for any $n,p<\omega$, we have $d(h_{n+p},h_n)=d(h_n^{-1}h_{n+p},1_G)=d(g_{n+1}^{-1}\cdots g_{n+p}^{-1},1_G)\le{\rm diam}(G_{n+1})\to 0$ as $n\to\infty$. It follows that $(h_n)$ is a $d$-Cauchy sequence in $G$, so it converges to some $h\in G$.

Let $x_\infty=h^{-1}\cdot x_0$. Since $h_n\to h$, we have $h_n^{-1}\to h^{-1}$, and hence $h_n^{-1}\cdot x_0\to h^{-1}\cdot x_0=x_\infty$. Note that $X$ is discrete, so there exists an integer $N$ such that $h_n^{-1}\cdot x_0=x_\infty$ for any $n>N$, thus $x_\infty=g_ng_{n-1}\cdots g_0\cdot x_0\in C_{n+1}$. This implies that $x_\infty\in\bigcap_nC_n$ and $G_n\cdot x_\infty=C_n$ for each $n<\omega$.

Finally, define $G_{x_\infty}=\{g\in G:g\cdot x_\infty=x_\infty\}$ and put $f:G\to X$ as $f(g)=g\cdot x_\infty$. Since $f$ is continuous and $\{x_\infty\}$ is clopen in $X$, it follows that $G_{x_\infty}=f^{-1}(x_\infty)$ is a clopen subgroup of $G$. So there exists an $m<\omega$ such that $G_m\subseteq G_{x_\infty}$. We now have $C_m=G_m\cdot x_\infty=\{x_\infty\}$, which is a singleton $G_m$-orbit, contradicting that $(m,C_m)\in T_\mathcal G^X$.
\end{proof}

Given two sets $X$ and $Y$. Let $\mathcal E=(E_n)$ and $\mathcal F=(F_n)$ be two decreasing sequences of equivalence relations on $X$ and $Y$ respectively. Let $\theta:X\to Y$ be an injection such that $\theta$ is a reduction of $E_n$ to $F_n$ for each $n<\omega$, i.e.,
$$\forall n<\omega\,\forall x,x'\in X\,(xE_nx'\iff\theta(x)F_n\theta(x')).$$
For $(n,C)\in T_\mathcal E^X$, define $\phi(n,C)=(n,[\theta(C)]_{F_n})$.

\begin{proposition}\label{embedding}
$\phi$ is a Lipschitz embedding from $T_\mathcal E^X$ to $T_\mathcal F^Y$. In particular, if $\theta$ is a bijection, then $\phi$ is an order-preserving isomorphism.
\end{proposition}

\begin{proof}
Note that $\theta$ is injective. Let us prove that $\phi(n,C)\in T_\mathcal{F}^Y$ holds for an arbitrary $(n,C)\in T_\mathcal{E}^X$. Indeed, it suffices to note that since $C$ is not a singleton, neither is $[\phi(C)]_{F_n}$. This proves $\phi(n,C)\in T_\mathcal{F}^Y$.
The rest of the proof is trivial.
\end{proof}

\begin{definition}
Let $(T,<)$ be a tree, $(n_i)$ a strictly increasing sequence of natural numbers. We define
$$T|(n_i)=\bigcup_iL_{n_i}(T).$$
Note that $(T|(n_i),<)$ is also a tree, and that $L_j(T|(n_i))=L_{n_j}(T)$ holds for each $j<\omega$. We call $T|(n_i)$ a {\it level-subtree} of $T$.
\end{definition}

\begin{lemma}
Let $(T,<)$ be a well-founded tree, $(n_i)$ a strictly increasing sequence of natural numbers. Then we have
$$\omega(\rho(T))\le\rho(T|(n_i))\le\rho(T).$$
In particular, if $\rho(T)$ is a limit ordinal, then $\rho(T|(n_i))=\rho(T)$.
\end{lemma}

\begin{proof}
$\rho(T|(n_i))\le\rho(T)$ follows from Proposition~\ref{subtree}. We prove $\omega(\rho(T))\le\rho(T|(n_i))$ by induction on $\rho(T)$.

First, if $\rho(T)<\omega$, then $\rho(T)=\min\{n:L_n(T)=\emptyset\}$, and hence $\rho(T|(n_i))=\min\{i:L_{n_i}(T)=\emptyset\}$. So we have $\omega(\rho(T))=0\le\rho(T|(n_i))$.

For $t\in T|(n_i)$, note that $(T|(n_i))_t=\{u\in T|(n_i):t=u\vee t<u\}$ is a level-subtree of $T_t$ as well.

{\sl Case 1:} If $\rho(T)$ is a limit ordinal, then $\omega(\rho(T))=\rho(T)$. Proposition~\ref{tree} implies that $\rho(T)=\sup\{\rho(T_t):t\in T\}$.
Since $\rho(T)$ is a limit ordinal and $\rho(T_t)$ is a successor ordinal, we have $\rho(T_t)<\rho(T)$ for $t\in T$.

{\sl Subcase 1.1:} If there is no maximum in $\{\omega(\rho(T_t)):t\in T\}$, we have
$$\rho(T)=\sup\{\rho(T_t):t\in T\}=\sup\{\omega(\rho(T_t)):t\in T\}.$$
By the inductive hypothesis, we have $\omega(\rho(T_t))\le\rho((T|(n_i))_t)$. Proposition~\ref{subtree} gives $\rho((T|(n_i))_t)\le\rho(T|(n_i))$ for each $t\in T$, so we have $\rho(T|(n_i))=\rho(T)$.

{\sl Subcase 1.2:} Otherwise, let $\alpha=\max\{\omega(\rho(T_t)):t\in T\}$. Since $\rho(T_t)<\rho(T)$ for $t\in T$, we have $\rho(T)=\alpha+\omega$. We can find a sequence $t_m,\,m<\omega$ in $L_0$ such that $\rho(T_{t_m})=\alpha+k_m$ with $\sup\{k_m:m<\omega\}=\omega$. By Proposition~\ref{tree}, for each $m<\omega$ and $n<k_m$ we can find $t_m^n\in L_n(T)$ such that $t_m=t_m^0<t_m^1<\dots<t_m^{k_m-1}$ and $\rho(T_{t_m^n})=\alpha+(k_m-n)$. For $k_m>n_0$, let $i_m$ be the largest $i$ such that $n_i<k_m$, then $t_m^{n_j}\in L_j(T|(n_i))=L_{n_j}(T)$ for $j\le i_m$. By the inductive hypothesis, $\alpha=\omega(\rho(T_{t_m^{n_j}}))\le\rho((T|(n_i))_{t_m^{n_j}})$. Since $\rho((T|(n_i))_{t_m^{n_j}})\ge\rho((T|(n_i))_{t_m^{n_{j+1}}})+1$ for each $j<i_m$, we have $\rho((T|(n_i))_{t_m^{n_0}})\ge\alpha+i_m$. By the definition of $i_m$, we have $\sup\{i_m:m<\omega\}=\omega$. This gives $\rho(T|(n_i))=\alpha+\omega=\rho(T)$.

{\sl Case 2:} If $\rho(T)=\omega(\rho(T))+n$ with $1\le n<\omega$, then there exists some $t_0\in L_0(T)$ such that $\rho(T)=\rho_T(t_0)+1$. Since $\{u\in T|(n_i):t_0<u\}$ is a level-subtree of $\{u\in T:t_0<u\}$ and $\rho(\{u\in T:t_0<u\})=\rho_T(t_0)<\rho(T)$, by the inductive hypothesis and Proposition~\ref{subtree}, we have
$$\omega(\rho_T(t_0))=\omega(\rho(\{u\in T:t_0<u\}))\le\rho(\{u\in T|(n_i):t_0<u\})\le\rho(T|(n_i)).$$
It follows that $\omega(\rho(T))=\omega(\rho_T(t_0))\le\rho(T|(n_i))$.
\end{proof}

In general, the tree $T_\mathcal G^X$ and the ordinal $\rho(T_\mathcal G^X)$ depend on $\mathcal G$, and not only on the action $G\curvearrowright X$. The following key lemma shows that $\omega(\rho(T_\mathcal G^X))$ is independent from the choice of $\mathcal G$.

\begin{lemma}\label{two dgnb}
Let $G$ be a non-archimedean CLI Polish group, $X$ be a countable discrete Polish $G$-space, and let $\mathcal G=(G_n),\mathcal G'=(G_n')\in{\rm dgnb}(G)$. Then $$\omega(\rho(T_\mathcal G^X))=\omega(\rho(T_{\mathcal G'}^X)).$$
\end{lemma}

\begin{proof}
(1) First, we consider the case where $(G_n')$ is a subsequence of $(G_n)$, i.e., there is a strictly increasing sequence $(n_i)$ of natural numbers such that $G_i'=G_{n_i}$ for each $i<\omega$.

We define $\psi:T_{\mathcal G'}^X\to T_\mathcal G^X$ as $\psi(i,C)=(n_i,C)$. It is clear that $\psi$ is an order-preserving isomorphism from $T_{\mathcal G'}^X$ onto $T_\mathcal G^X|(n_i)$. It follows that
$$\omega(\rho(T_\mathcal G^X))\le\rho(T_{\mathcal G'}^X)=\rho(T_\mathcal G^X|(n_i))\le\rho(T_\mathcal G^X).$$
So we have $\omega(\rho(T_\mathcal G^X))=\omega(\rho(T_{\mathcal G'}^X))$.

(2) Since $(G_n),(G_n')\in{\rm dgnb}(G)$, we can find two strictly increasing natural numbers $(n_i)$ and $(m_j)$ such that $n_0=0,m_0=0$, and
$$G_0\supseteq G_{m_0}'\supseteq G_{n_1}\supseteq G_{m_1}'\supseteq G_{n_2}\supseteq\cdots.$$
Define $H_{2i}=G_{n_i}$ and $H_{2i+1}=G_{m_i}'$ for each $i<\omega$. Then $(H_k)\in{\rm dgnb}(G)$. Put $\mathcal H=(H_k),\mathcal K=(G_{n_i})$, and $\mathcal K'=(G_{m_i}')$.

Note that $(G_{n_i})$ is a subsequence of $(G_n)$ and also a subsequence of $(H_k)$. From (1), we obtain
$$\omega(\rho(T_\mathcal G^X))=\omega(\rho(T_\mathcal K^X))=\omega(\rho(T_\mathcal H^X)).$$
Similarly, we obtain
$$\omega(\rho(T_{\mathcal G'}^X))=\omega(\rho(T_{\mathcal K'}^X))=\omega(\rho(T_\mathcal H^X)).$$
So we have $\omega(\rho(T_\mathcal G^X))=\omega(\rho(T_{\mathcal G'}^X))$.
\end{proof}

Now we are going to find a special $G$-space $X(\mathcal G)$ such that $\omega(\rho(T_\mathcal G^{X(\mathcal G)}))$ attains the maximum value among all $\omega(\rho(T_\mathcal G^X))$. This leads to the conclusion that the value of $\omega(\rho(T_\mathcal G^{X(\mathcal G)}))$ is determined by $G$ itself.

\begin{lemma}\label{surjection-tree}
Given two sets $X$ and $Y$. Let $\mathcal E=(E_n)$ and $\mathcal F=(F_n)$ be two decreasing sequences of equivalence relations on $X$ and $Y$ respectively. Let $\theta:X\to Y$ be a surjection such that
$$\forall n<\omega\,\forall x\in X\,(\theta([x]_{E_n})=[\theta(x)]_{F_n}).$$
Then there exists a Lipschitz embedding $\psi:T_\mathcal F^Y\to T_\mathcal E^X$. In particular, if $T_\mathcal E^X$ is well-founded, so is $T_\mathcal F^Y$, and $\rho(T_\mathcal F^Y)\le\rho(T_\mathcal E^X)$ holds.
\end{lemma}

\begin{proof}
For any $(n,C)\in T_\mathcal F^Y$, we construct $\psi(n,C)$ by induction on $n$ such that $\psi(n,C)=(n,[x]_{E_n})$ for some $x\in X$ with $[\theta(x)]_{F_n}=C$.

For $n=0$, since $\theta$ is a surjection, we can find an $x\in X$ with $\theta(x)\in C$. Then we let $\psi(0,C)=(0,[x]_{E_0})$.

For $n>0$, since $C$ is an $F_n$-equivalence class, there exists a unique $F_{n-1}$-equivalence class $D\supseteq C$. By the inductive hypothesis, we can find an $x'\in X$ such that $\psi(n-1,D)=(n-1,[x']_{E_{n-1}})$ with $[\theta(x')]_{F_{n-1}}=D$. Since $\theta([x']_{E_{n-1}})=[\theta(x')]_{F_{n-1}}=D\supseteq C$, we can find an $x\in[x']_{E_{n-1}}$ such that $\theta(x)\in C$. Then we put $\psi(n,C)=(n,[x]_{E_n})$.

Since $C$ is not a singleton, the equality $\theta([x]_{E_n})=[\theta(x)]_{F_n}=C$ implies that $[x]_{E_n}$ is not a singleton. So $\psi(n,C)\in T_\mathcal E^X$. From the construction, it is routine to check that $\psi:T_\mathcal F^Y\to T_\mathcal E^X$ is a Lipschitz embedding.

Finally, if $T_\mathcal E^X$ is well-founded, then by Proposition~\ref{subtree}, $T_\mathcal F^Y$ is well-founded as well, and $\rho(T_\mathcal F^Y)\le\rho(T_\mathcal E^X)$.
\end{proof}

\begin{definition}
Let $G$ be a non-archimedean CLI Polish group, and let $\mathcal G=(G_n)\in{\rm dgnb}(G)$. For each $k<\omega$, we define an action $G\curvearrowright G/G_k$ as $g\cdot hG_k=ghG_k$ for $g,h\in G$, and let $\rho^k(\mathcal G)$ denote $\rho(T^{G/G_k}_\mathcal G)$. Afterwards, we let $X(\mathcal G)=\bigcup_kG/G_k$ and $\rho(\mathcal G)=\rho(T^{X(\mathcal G)}_\mathcal G)$.
\end{definition}

Note that $G\cdot gG_k=\{hgG_k:h\in G\}=G/G_k$ for any $g\in G$.

It is clear that $\rho^0(\mathcal G)=0$, since $G=G_0$.

\begin{lemma}\label{rho}
\begin{enumerate}
\item[(1)] $\rho(\mathcal G)=\sup\{\rho^k(\mathcal G):k<\omega\}$.
\item[(2)] $(\rho^k(\mathcal G))$ is an increasing sequence of countable non-limit ordinals.
\end{enumerate}
\end{lemma}

\begin{proof}
(1) Note that $L_0(T^{X(\mathcal G)}_\mathcal G)=\{(0,G/G_k):G\ne G_k\}$. We point out that 
$(T^{X(\mathcal G)}_\mathcal G)_{(0,G/G_k)}\cong T^{G/G_k}_\mathcal G$
for $G\ne G_k$, and $(T^{X(\mathcal G)}_\mathcal G)_{(0,G/G_k)}=T^{G/G_k}_\mathcal G=\emptyset$ for $G=G_k$. So
$$\begin{array}{ll}\rho(\mathcal G)&=\rho(T^{X(\mathcal G)}_\mathcal G)=\sup\{\rho((T^{X(\mathcal G)}_\mathcal G)_{(0,G/G_k)}):k<\omega\}\cr
&=\sup\{\rho(T^{G/G_k}_\mathcal G):k<\omega\}=\sup\{\rho^k(\mathcal G):k<\omega\}.\end{array}$$

(2) Given $k<\omega$, we define $\theta:G/G_{k+1}\to G/G_k$ as $\theta(gG_{k+1})=gG_k$ for $g\in G$. It is clear that $\theta$ is well-defined and is a surjection. Furthermore, for $n<\omega$ and $g\in G$, we have
$$\begin{array}{ll}\theta(G_n\cdot gG_{k+1})&=\{\theta(hgG_{k+1}):h\in G_n\}\cr
&=\{hgG_k:h\in G_n\}=G_n\cdot gG_k=G_n\cdot\theta(gG_{k+1}).\end{array}$$
From Lemma~\ref{surjection-tree}, we have
$$\rho(T_\mathcal G^{G/G_k})\le\rho(T_\mathcal G^{G/G_{k+1}}),$$
i.e., $(\rho^k(\mathcal G))$ is increasing.

For each $k<\omega$, since $T_\mathcal G^{G/G_k}$ is countable, $\rho^k(\mathcal G)$ is countable as well. If $G=G_k$, then $T^{G/G_k}_\mathcal G=\emptyset$; otherwise if $G\ne G_k$, then $L_0(T^{G/G_k}_\mathcal G)=\{(0,G/G_k)\}$ is a singleton. So $\rho^k(\mathcal G)=\rho(T^{G/G_k}_\mathcal G)$ is either $0$ or a successor ordinal.
\end{proof}

Recall that a $G$-space $X$ is said to be \emph{transitive} if $X$ itself is an orbit.

\begin{lemma}\label{transitive}
Let $G$ be a non-archimedean CLI Polish group, $X$ a countable discrete transitive Polish $G$-space, and let $\mathcal G=(G_n)\in{\rm dgnb}(G)$. Then we can find some $k<\omega$ such that $\rho(T_\mathcal G^X)\le\rho^k(\mathcal G)$.
\end{lemma}

\begin{proof}
Fix an $x\in X$. Since $\{x\}$ is clopen in $X$, by the continuity of the group action of $G$ on $X$, we know $G_x$ is a clopen subgroup of $G$. So there is some $k<\omega$ such that $G_k\subseteq G_x$. Then we can define $\theta:G/G_k\to X$ as $\theta(gG_k)=g\cdot x$ for $g\in G$. Since $X$ is a transitive $G$-space, $\theta$ is surjective and $\theta(G_n\cdot gG_k)=G_n\cdot\theta(gG_k)$ for each $n<\omega$ and $g\in G$. By Lemma~\ref{surjection-tree}, we have $\rho(T_\mathcal G^X)\le\rho^k(\mathcal G)$.
\end{proof}

\begin{corollary}\label{non-transitive}
Let $G$ be a non-archimedean CLI Polish group, $X$ a countable discrete Polish $G$-space, and let $\mathcal G=(G_n)\in{\rm dgnb}(G)$. Then we have $$\rho(T_\mathcal G^X)\le\rho(\mathcal G).$$
\end{corollary}

\begin{proof}
Note that $L_0(T_\mathcal G^X)=\{(0,G\cdot x):x\in X\wedge G\cdot x\ne\{x\}\}$ and $T_\mathcal G^{G\cdot x}\cong(T_\mathcal G^X)_{(0,G\cdot x)}$ for $G\cdot x\ne\{x\}$, so we have $$\rho(T_\mathcal G^X)=\sup\{\rho((T_\mathcal G^X)_{(0,G\cdot x)}):x\in X\}=\sup\{\rho(T_\mathcal G^{G\cdot x}):x\in X\}.$$
By Lemma~\ref{transitive}, we have
$$\rho(T_\mathcal G^X)\le\sup\{\rho^k(\mathcal G):k<\omega\}=\rho(\mathcal G).$$
\end{proof}

\begin{corollary}
Let $G$ be a non-archimedean Polish group, $\mathcal G=(G_n)\in{\rm dgnb}(G)$. Then $G$ is CLI iff $T_\mathcal G^{G/G_k}$ is well-founded for any $k<\omega$.
\end{corollary}

\begin{proof}
The $(\Rightarrow)$ part follows from Lemma~\ref{CLItoWF}.

$(\Leftarrow)$. Given a countable Polish $G$-space $X$. Following the arguments in the proof of Lemma~\ref{transitive}, we can see that, for any $x\in X$, there is a $k<\omega$ and a Lipschitz embedding from $T^{G\cdot x}_\mathcal G$ to $T_\mathcal G^{G/G_k}$. Since $T_\mathcal G^{G/G_k}$ is well-founded, so is $T^{G\cdot x}_\mathcal G$. By the arbitrariness of $x\in X$, we conclude that $T^X_\mathcal G$ is also well-founded. Combining this with \cite[Theorem 6]{malicki11}, it follows that $G$ is CLI.
\end{proof}

\begin{theorem}\label{omega=omega}
Let $G$ be a non-archimedean CLI Polish group, and $\mathcal G=(G_n),\mathcal G'=(G_n')\in{\rm dgnb}(G)$. Then $\omega(\rho(\mathcal G))=\omega(\rho(\mathcal G'))$ holds.
\end{theorem}

\begin{proof}
By Corollary~\ref{non-transitive}, we have $\rho(T^{X(\mathcal G)}_{\mathcal G'})\le\rho(\mathcal G')$. Thus Lemma~\ref{two dgnb} implies
$$\omega(\rho(\mathcal G))=\omega(\rho(T^{X(\mathcal G)}_\mathcal G))=\omega(\rho(T^{X(\mathcal G)}_{\mathcal G'}))\le\omega(\rho(\mathcal G')),$$
and vice versa.
\end{proof}

By the preceding theorem, there is a unique ordinal $\beta<\omega_1$, which is independent from the choice of $\mathcal G=(G_n)\in{\rm dgnb}(G)$, such that
$$\omega(\rho(\mathcal G))=\omega\cdot\beta.$$
Consequently, we can define a rank for a given non-archimedean CLI Polish group as follows:

\begin{definition}
    Suppose $G$ is a non-archimedean CLI Polish group. We define 
    $$\mathrm{rank}(G)=\textrm{the unique }\beta \textrm{ such that } \omega(\rho(\mathcal{G}))=\omega\cdot\beta$$
    for some $\mathcal{G}\in\mathrm{dgnb}(G)$.
\end{definition}

\begin{lemma}\label{rank}
\begin{enumerate}
\item[(1)] If $\rho(\mathcal G)=\omega\cdot{\rm rank}(G)$, then either ${\rm rank}(G)=0$ or $\rho^k(\mathcal G)<\omega\cdot{\rm rank}(G)$ for any $k<\omega$.
\item[(2)] If $\rho(\mathcal G)>\omega\cdot{\rm rank}(G)$, then there exists an $m>0$ such that $\rho^k(\mathcal G)=\omega\cdot{\rm rank}(G)+m$ for large enough $k<\omega$.
\end{enumerate}
\end{lemma}

\begin{proof}
(1) If ${\rm rank}(G)>0$, then $\omega\cdot{\rm rank}(G)$ is a limit ordinal. By Lemma~\ref{rho}, $\rho^k(\mathcal G)$ is either $0$ or a successor ordinal, so $\rho^k(\mathcal G)<\omega\cdot{\rm rank}(G)$ for any $k<\omega$.

(2) Clearly, $\rho(\mathcal G)=\omega\cdot{\rm rank}(G)+m$ for some $0<m<\omega$. Again by Lemma~\ref{rho}, $(\rho^k(\mathcal G))$ is increasing, so $\rho^k(\mathcal G)=\omega\cdot{\rm rank}(G)+m$ for large enough $k<\omega$.
\end{proof}

\begin{theorem}\label{equi-omega}
Let $G$ be a non-archimedean CLI Polish group, $\mathcal G=(G_n),\mathcal G'=(G_n')\in{\rm dgnb}(G)$. Then $\rho(\mathcal G)=\omega\cdot{\rm rank}(G)$ iff $\rho(\mathcal G')=\omega\cdot{\rm rank}(G)$.
\end{theorem}

\begin{proof}
If ${\rm rank}(G)=0$, then $\rho(\mathcal G)=\omega\cdot{\rm rank}(G)$ implies that $\rho^k(\mathcal G)=0$ for any $k<\omega$. So $T_\mathcal G^{G/G_k}=\emptyset$, and hence $G_0\cdot G_k\notin T_\mathcal G^{G/G_k}$. It follows that $G_0\cdot G_k=\{G_k\}$, i.e., $G=G_0=G_k$ for any $k<\omega$. This implies that $G=\{1_G\}$. Then we can easily see that $T_{\mathcal G'}^{G/G_k'}=\emptyset$ for $k<\omega$. So $\rho(\mathcal G')=\omega\cdot{\rm rank}(G)$ holds. And vice versa.

If ${\rm rank}(G)>0$, assume for contradiction that $\rho(\mathcal G)=\omega\cdot{\rm rank}(G)$, but $\rho(\mathcal G')>\omega\cdot{\rm rank}(G)$. From Lemma~\ref{rank}, we have $\rho^k(\mathcal G)<\omega\cdot{\rm rank}(G)$ for any $k<\omega$, but $\rho^l(\mathcal G')=\omega\cdot{\rm rank}(G)+m$ for some $0<m<\omega$ and large enough $l<\omega$. From Lemma~\ref{two dgnb} and Lemma~\ref{transitive}, we can conclude that for any $l<\omega$, there is some $k<\omega$ such that $$\omega(\rho^l(\mathcal G'))=\omega(\rho(T_{\mathcal G'}^{G/G_l'}))=\omega(\rho(T_\mathcal G^{G/G_l'}))\le\rho(T_\mathcal G^{G/G_l'})\le\rho^k(\mathcal G).$$
A contradiction!
\end{proof}

Now we are ready to define a hierarchy on the class of non-archimedean CLI Polish groups.

\begin{definition}
Let $G$ be a non-archimedean CLI Polish group, $\mathcal G=(G_n)\in{\rm dgnb}(G)$, and let $\alpha<\omega_1$ be an ordinal.
\begin{enumerate}
\item[(1)] If $\rho(\mathcal G)\le\omega\cdot\alpha$, we say that $G$ is {\it $\alpha$-CLI};
\item[(2)] if $\omega(\rho(\mathcal G))\le\omega\cdot\alpha$, i.e., ${\rm rank}(G)\le\alpha$, we say that $G$ is {\it L-$\alpha$-CLI}.
\end{enumerate}
It is clear that, if $G$ is L-$\alpha$-CLI, then it is also $(\alpha+1)$-CLI.

From Theorem~\ref{omega=omega} and Theorem~\ref{equi-omega}, we see that the definitions of $\alpha$-CLI and L-$\alpha$-CLI are independent from the choice of $\mathcal G\in{\rm dgnb}(G)$.
\end{definition}

Recall that a metric $d$ on a group $G$ is {\it two-sided invariant} if $d(gh,gk)=d(h,k)=d(hg,kg)$ for all $g,h,k\in G$. A Polish group is {\it TSI} if it admits a compatible complete two-sided invariant metric.

\begin{theorem}\label{0-1}
Let $G$ be a non-archimedean CLI Polish group. The following hold:
\begin{enumerate}
\item[(1)] $G$ is $0$-CLI iff $G=\{1_G\}$;
\item[(2)] $G$ is L-$0$-CLI iff $G$ is discrete; and
\item[(3)] $G$ is $1$-CLI iff $G$ is TSI.
\end{enumerate}
\end{theorem}

\begin{proof}
Fix a sequence $\mathcal G=(G_n)\in{\rm dgnb}(G)$.

(1) It follows from the first paragraph of the proof of Theorem~\ref{equi-omega}.

(2) If $G$ is L-$0$-CLI, then ${\rm rank}(G)=0$. So there is an $m<\omega$ such that $\rho(T_\mathcal G^{G/G_k})=\rho^k(\mathcal G)=m$ for large enough $k<\omega$. Then we have $L_m(T_\mathcal G^{G/G_k})=\emptyset$. This implies that $G_m\cdot G_k=\{G_k\}$. So $G_m\subseteq G_k$ for large enough $k<\omega$, and thus $G_m=\{1_G\}$. It follows that $G$ is discrete.

On the other hand, if $G$ is discrete, then there is an $m<\omega$ such that $G_m=\{1_G\}$. Therefore, for any $k<\omega$, we have $L_m(T_\mathcal G^{G/G_k})=\emptyset$, and hence $\rho^k(\mathcal G)\le m$. This implies that ${\rm rank}(G)=0$, i.e., $G$ is L-$0$-CLI.

(3) If $G$ is $1$-CLI, then Lemma~\ref{rank} implies that $\rho^k(\mathcal G)<\omega$ for any $k<\omega$. So there is an $m_k<\omega$ such that $L_{m_k}(T_\mathcal G^{G/G_k})=\emptyset$. It follows that, for $g\in G$, $G_{m_k}\cdot gG_k=\{gG_k\}$, so $g^{-1}G_{m_k}g\subseteq G_k$. Put $U_k=\bigcup\{g^{-1}G_{m_k}g:g\in G\}\subseteq G_k$. We can see that $(U_k)$ is a neighborhood basis of $1_G$ with $g^{-1}U_kg=U_k$ for all $g\in G$. By Klee's theorem (c.f.~\cite{klee} or~\cite[Exercise 2.1.4]{gaobook}), $G$ is TSI.

On the other hand, if $G$ is TSI, again by Klee's theorem, we can find a neighborhood basis $(U_m)$ of $1_G$ with $g^{-1}U_mg=U_m$ for all $g\in G$. For any $n<\omega$, there is an $m_n<\omega$ such that $U_{m_n}\subseteq G_n$. Let $V_n=U_{m_n}\cap U_{m_n}^{-1}$ and $G_n'=\bigcup_iV_n^i$. Then $G_n'$ is an open normal subgroup of $G$ with $G_n'\subseteq G_n$. So $(G_n')\in{\rm dgnb}(G)$. Let $\mathcal G'=(G_n')$. Then $G_k'\cdot gG_k'=\{gG_k'\}$ for all $g\in G$ and $k<\omega$, thus $L_k(T_{\mathcal G'}^{G/G_k'})=\emptyset$. So $\rho^k(\mathcal G')\le k<\omega$, and hence $G$ is $1$-CLI.
\end{proof}

Clause (2) in the preceding theorem can be generalized to all $\alpha<\omega_1$.

\begin{definition}
Let $G$ be a non-archimedean CLI Polish group, and $\alpha<\omega_1$. We say that $G$ is {\it locally $\alpha$-CLI} if $G$ has an open subgroup which is $\alpha$-CLI.
\end{definition}

\begin{theorem}\label{locally}
Let $G$ be a non-archimedean CLI Polish group, and $\alpha<\omega_1$. Then $G$ is L-$\alpha$-CLI iff $G$ is locally $\alpha$-CLI.
\end{theorem}

\begin{proof}
$(\Rightarrow)$. Suppose $G$ is L-$\alpha$-CLI. Without loss of generality, we may assume that $G$ is not $\alpha$-CLI. Fix a sequence $\mathcal G=(G_n)\in{\rm dgnb}(G)$. There exists an $m\ge 1$ such that $\rho(\mathcal G)=\omega\cdot\alpha+m$, and thus we can pick a $k_0>m$ such that $\rho^k(\mathcal G)=\omega\cdot\alpha+m$ for any $k\ge k_0$. We will show that $G_{k_0}$ is $\alpha$-CLI.

Let $H=G_{k_0}$ and $H_n=G_{n+k_0}$ for $n<\omega$. Then $(H_n)\in{\rm dgnb}(H)$. Put $\mathcal H=(H_n)$. Given $k<\omega$, define $\phi:T^{H/H_k}_\mathcal H\to(T^{G/G_{k+k_0}}_\mathcal G)_{(k_0,G_{k_0}/G_{k+k_0})}$ as $\phi(n,C)=(n+k_0,C)$. It is trivial to see that $\phi$ is an order-preserving isomorphism. From Proposition~\ref{L_k(T)}, since $k_0>m$, we have
$$\rho^k(\mathcal H)=\rho(T^{H/H_k}_\mathcal H)=\rho((T^{G/G_{k+k_0}}_\mathcal G)_{(k_0,G_{k_0}/G_{k+k_0})})\le\omega\cdot\alpha.$$
So $\rho(\mathcal H)\le\omega\cdot\alpha$, and hence $H=G_{k_0}$ is $\alpha$-CLI.

$(\Leftarrow)$. Suppose $G$ is locally $\alpha$-CLI. Let $H$ be an open subgroup of $G$ which is $\alpha$-CLI, and let $\mathcal H=(H_n)\in{\rm dgnb}(H)$. Then $\rho^k(\mathcal H)\le\omega\cdot\alpha$ for $k<\omega$.

Put $G_0=G$ and $G_n=H_{n-1}$ for $n\ge 1$. Then $(G_n)\in{\rm dgnb}(G)$. Put $\mathcal G=(G_n)$. Given $g\in G$ and $k<\omega$, by similar arguments in the $(\Rightarrow)$ part, we have $T^{H\cdot gH_k}_\mathcal H\cong(T^{G/G_{k+1}}_\mathcal G)_{(1,G_1\cdot gG_{k+1})}$. By Lemma~\ref{transitive}, there exists an $l<\omega$ such that $\rho(T^{H\cdot gH_k}_\mathcal H)\le\rho^l(\mathcal H)$. Therefore, by Proposition~\ref{tree},
$$\begin{array}{ll}\rho^{k+1}(\mathcal G)&=\rho(T^{G/G_{k+1}}_\mathcal G)\cr
&\le\sup\{\rho((T^{G/G_{k+1}}_\mathcal G)_{(1,G_1\cdot gG_{k+1})}):g\in G\}+1\cr
&=\sup\{\rho(T^{H\cdot gH_k}_\mathcal H):g\in G\}+1\cr
&\le\sup\{\rho^l(\mathcal H):l<\omega\}+1\cr
&\le\omega\cdot\alpha+1.\end{array}$$
So $\rho(\mathcal G)\le\omega\cdot\alpha+1$, and hence $G$ is L-$\alpha$-CLI.
\end{proof}

\section{Properties of the hierarchy}

\begin{theorem}\label{subgroup}
Let $G$ be a non-archimedean CLI Polish group, $H$ a closed subgroup of $G$, and $\alpha<\omega_1$. If $G$ is $\alpha$-CLI (or L-$\alpha$-CLI), so is $H$. In particular, we have ${\rm rank}(H)\le{\rm rank}(G)$.
\end{theorem}

\begin{proof}
Let $\mathcal G=(G_n)\in{\rm dgnb}(G)$, and put $H_n=H\cap G_n$ for $n<\omega$. It is clear that $(H_n)\in{\rm dgnb}(H)$. Put $\mathcal H=(H_n)$. We only need to show that $\rho(\mathcal H)\le\rho(\mathcal G)$.

Given $k<\omega$, define $\theta:H/H_k\to G/G_k$ as $\theta(hH_k)=hG_k$ for $h\in H$. By Proposition~\ref{embedding}, there is a Lipschitz embedding from $T^{H/H_k}_\mathcal H$ to $T^{G/G_k}_\mathcal G$. So
$$\rho^k(\mathcal H)=\rho(T^{H/H_k}_\mathcal H)\le\rho(T^{G/G_k}_\mathcal G)=\rho^k(\mathcal G).$$
Then we have $\rho(\mathcal H)\le\rho(\mathcal G)$ as desired.
\end{proof}

\begin{theorem}
Let $G$ be a non-archimedean CLI Polish group, $N$ a closed normal subgroup of $G$, and $\alpha<\omega_1$. If $G$ is $\alpha$-CLI (or L-$\alpha$-CLI), so is $G/N$. In particular, we have ${\rm rank}(G/N)\le{\rm rank}(G)$.
\end{theorem}

\begin{proof}
Let $\mathcal G=(G_n)\in{\rm dgnb}(G)$, and put $H_n=G_n\cdot N=\{\hat gN:\hat g\in G_n\}$ for $n<\omega$. It is clear that $H_0=G/N$ and $(H_n)\in{\rm dgnb}(G/N)$. Put $\mathcal H=(H_n)$. We only need to show that $\rho(\mathcal H)\le\rho(\mathcal G)$.

Given $k<\omega$, define $\theta:G/G_k\to(G/N)/H_k$ as $\theta(gG_k)=(gN)H_k$ for $g\in G$. Note that for $n<\omega$,
$$\begin{array}{ll}\theta(G_n\cdot gG_k)&=\theta(\{\hat ggG_k:\hat g\in G_n\})\cr
&=\{(\hat ggN)H_k:\hat g\in G_n\}\cr
&=\{(\hat gN)(gN)H_k:\hat g\in G_n\}\cr
&=\{(\hat gN)\theta(gG_k):\hat g\in G_n\}\cr
&=\{\hat gN:\hat g\in G_n\}\theta(gG_k)=H_n\cdot\theta(gG_k).\end{array}$$
By Lemma~\ref{surjection-tree}, we have
$$\rho^k(\mathcal H)=\rho(T^{(G/N)/H_k}_\mathcal H)\le\rho(T^{G/G_k}_\mathcal G)=\rho^k(\mathcal G).$$
So $\rho(\mathcal H)\le\rho(\mathcal G)$ holds as desired.
\end{proof}

The above two theorems involve closed subgroups and quotient groups. Now we turn to product groups, which are more complicated. We discuss finite product groups first.

\begin{lemma}\label{XtimesY}
Let $X,Y$ be two sets, $\mathcal E=(E_n)$ and $\mathcal F=(F_n)$ be two decreasing sequences of equivalence relations on $X$ and $Y$ respectively. Let $\mathcal{E\times F}=(E_n\times F_n)$.
\begin{enumerate}
\item[(1)] $T_\mathcal{E\times F}^{X\times Y}$ is well-founded iff $T_\mathcal E^X$ and $T_\mathcal F^Y$ are well-founded.
\item[(2)] If $T_\mathcal{E\times F}^{X\times Y}$ is well-founded, then we have
$$\rho(T_\mathcal{E\times F}^{X\times Y})=\max\{\rho(T_\mathcal E^X),\rho(T_\mathcal F^Y)\}.$$
\end{enumerate}
\end{lemma}

\begin{proof}
First, note that $[(x,y)]_{E_n\times F_n}=[x]_{E_n}\times[y]_{F_n}$ for all $(x,y)\in X\times Y$ and $n<\omega$.

(1) For any sequences $(x_n),(y_n)$ in $X,Y$ respectively, $((n,[(x_n,y_n)]_{E_n\times F_n}))$ is an infinite branch of $T_\mathcal{E\times F}^{X\times Y}$ iff either $((n,[x_n]_{E_n}))$ or $((n,[y_n]_{F_n}))$ is an infinite branch of $T_\mathcal E^X$ or $T_\mathcal F^Y$ respectively.

(2) If $T_\mathcal{E\times F}^{X\times Y}$ is well-founded, by (1), $T_\mathcal E^X$ and $T_\mathcal F^Y$ are also well-founded. For all $(x,y)\in X\times Y$ and $n<\omega$, note that
$$\begin{array}{ll} & (n,[(x,y)]_{E_n\times F_n})\in T^{X\times Y}_\mathcal{E\times F}\cr
\iff & [(x,y)]_{E_n\times F_n}\ne\{(x,y)\}\cr
\iff & [x]_{E_n}\ne\{x\}\vee[y]_{F_n}\ne\{y\}\cr
\iff & (n,[x]_{E_n})\in T_\mathcal E^X\vee(n,[y]_{F_n})\in T_\mathcal F^Y.\end{array}$$
By Proposition~\ref{tree}, it is routine to prove
$$\rho((T_\mathcal{E\times F}^{X\times Y})_{(n,[(x,y)]_{E_n\times F_n})})=\max\{\rho((T_\mathcal E^X)_{(n,[x]_{E_n})}),\rho((T_\mathcal F^Y)_{(n,[y]_{F_n})})\}$$
by induction on $\rho((T_\mathcal{E\times F}^{X\times Y})_{(n,[(x,y)]_{E_n\times F_n})})$. Taking supremum on both sides of the above formula, we get
$$\rho(T_\mathcal{E\times F}^{X\times Y})=\max\{\rho(T_\mathcal E^X),\rho(T_\mathcal F^Y)\}.$$
\end{proof}

\begin{corollary}\label{GtimesH}
Let $G$ and $H$ be two non-archimedean CLI Polish groups, $\mathcal G=(G_n)\in{\rm dgnb}(G)$ and $\mathcal H=(H_n)\in{\rm dgnb}(H)$. Then we have $\mathcal{G\times H}=(G_n\times H_n)\in{\rm dgnb}(G\times H)$ and
$$\rho^k(\mathcal{G\times H})=\max\{\rho^k(\mathcal G),\rho^k(\mathcal H)\}\quad(\forall k<\omega),$$
$$\rho(\mathcal{G\times H})=\max\{\rho(\mathcal G),\rho(\mathcal H)\},$$
$${\rm rank}(G\times H)=\max\{{\rm rank}(G),{\rm rank}(H)\}.$$
\end{corollary}

\begin{proof}
For any $k<\omega$, put $X=G/G_k,Y=H/H_k$, and define $E_n,F_n$ on $X$ and $Y$ for each $n<\omega$ respectively by
$$(\hat gG_k,\tilde gG_k)\in E_n\iff\exists g\in G_n\,(g\hat gG_k=\tilde gG_k)\quad(\forall\hat g,\tilde g\in G),$$
$$(\hat hH_k,\tilde hH_k)\in F_n\iff\exists h\in H_n\,(h\hat hH_k=\tilde hH_k)\quad(\forall\hat h,\tilde h\in H).$$
Then Lemma~\ref{XtimesY} gives $\rho^k(\mathcal{G\times H})=\max\{\rho^k(\mathcal G),\rho^k(\mathcal H)\}$. The rest follows trivially.
\end{proof}

Now we are ready to discuss countably infinite product groups.

\begin{lemma}\label{times}
Let $(G^i)$ be a sequence of non-archimedean CLI Polish groups and $\mathcal G^i=(G^i_n)\in{\rm dgnb}(G^i)$ for each $i<\omega$. We define $G=\prod_iG^i$ and
$$G_n=\prod_{i<n}G^i_n\times\prod_{i\ge n}G^i\quad(\forall n<\omega).$$
Then $\mathcal G=(G_n)\in{\rm dgnb}(G)$ and for $k<\omega$, we have
$$\rho^k(\mathcal G)\le\max\{\rho^k(\mathcal G^i):i<k\}+k.$$
\end{lemma}

\begin{proof}
Given $k<\omega$, we put $Y=G^0/G^0_k\times\dots\times G^{k-1}/G^{k-1}_k$ and define $\theta:G/G_k\to Y$ as
$$\theta((g_i)G_k)=(g_0G^0_k,\dots,g_{k-1}G^{k-1}_k)$$
for $(g_i)\in G=\prod_iG^i$. By the definition of $G_k$, it is easy to see that $\theta$ is a bijection. Moreover, we put
$$H_n=\left\{\begin{array}{ll}\prod_{i<n}G^i_n\times\prod_{n\le i<k}G^i, & n<k,\cr
\prod_{i<k}G^i_n, & n\ge k,\end{array}\right.$$
then $\theta$ is a reduction of $E^{G/G_k}_{G_n}$ to $E^Y_{H_n}$ for each $n<\omega$. Put $\mathcal H=(H_n)$. By Proposition~\ref{embedding}, $(n,G_n\cdot (g_i)G_k)\mapsto(n,H_n\cdot\theta((g_i)G_k))$ is an order-preserving isomorphism from $T^{G/G_k}_\mathcal G$ to $T^Y_\mathcal H$. So $\rho^k(\mathcal G)=\rho(T^{G/G_k}_\mathcal G)=\rho(T^Y_\mathcal H)$.

For $n\ge k$ and $(g_i)\in G$, we have
$$H_n\cdot\theta((g_i)G_k)=(G^0_n\cdot g_0G^0_k)\times\dots\times(G^{k-1}_n\cdot g_{k-1}G^{k-1}_k).$$
Lemma~\ref{XtimesY} gives
$$\rho((T^Y_\mathcal H)_{(k,H_k\cdot\theta((g_i)G_k))})=\max\{\rho((T^{G^i/G^i_k}_{\mathcal G^i})_{(k,G^i_k\cdot g_iG^i_k)}):i<k\}.$$
Hence, by Proposition~\ref{L_k(T)}.(1),
$$\begin{array}{ll}\rho^k(\mathcal G)&=\rho(T^{G/G_k}_\mathcal G)=\rho(T^Y_\mathcal H)\cr
&\le\sup\{\rho((T^Y_\mathcal H)_{(k,H_k\cdot\theta((g_i)G_k))}):(g_i)\in G\}+k\cr
&=\sup\{\max\{\rho((T^{G^i/G^i_k}_{\mathcal G^i})_{(k,G^i_k\cdot g_iG^i_k)}):i<k\}:(g_i)\in G\}+k\cr
&=\max\{\sup\{\rho((T^{G^i/G^i_k}_{\mathcal G^i})_{(k,G^i_k\cdot gG^i_k)}):g\in G^i\}:i<k\}+k\cr
&\le\max\{\rho(T^{G^i/G^i_k}_{\mathcal G^i}):i<k\}+k\cr
&=\max\{\rho^k(\mathcal G^i):i<k\}+k.\end{array}$$
\end{proof}

\begin{theorem}\label{prod}
Let $(G^i)$ be a sequence of non-archimedean CLI Polish groups, $\alpha<\omega_1$, and let $G=\prod_iG^i$. Then we have
\begin{enumerate}
\item[(1)] $G$ is $\alpha$-CLI iff all $G^i$ are $\alpha$-CLI; and
\item[(2)] $G$ is L-$\alpha$-CLI iff all $G^i$ are L-$\alpha$-CLI and for all but finitely many $i$, $G^i$ is $\alpha$-CLI.
\end{enumerate}
\end{theorem}

\begin{proof}
Fix a $\mathcal G^i=(G^i_n)\in{\rm dgnb}(G^i)$ for each $i<\omega$. Put
$$G_n=\prod_{i<n}G^i_n\times\prod_{i\ge n}G^i\quad(\forall n<\omega).$$

(1) Assume $G$ is $\alpha$-CLI. Since each $G^i$ is topologically isomorphic to a closed subgroup of $G$, by Theorem~\ref{subgroup}, $G^i$ is also $\alpha$-CLI.

On the other hand, assume all $G^i$ are $\alpha$-CLI. Now by Lemma~\ref{rank}.(1), $\rho^k(\mathcal G^i)<\omega\cdot\alpha$ holds for all $i,k<\omega$. Then Lemma~\ref{times} implies
$$\rho^k(\mathcal G)\le\max\{\rho^k(\mathcal G^i):i<k\}+k<\omega\cdot\alpha$$
holds for each $k<\omega$. Consequently, $G$ is $\alpha$-CLI.

(2) Assume $G$ is L-$\alpha$-CLI. Since each $G^i$ is topologically isomorphic to a closed subgroup of $G$, we can see $G^i$ is L-$\alpha$-CLI too. Moreover, by Theorem~\ref{locally}, there is an open subgroup $H$ of $G$ which is $\alpha$-CLI. Hence, $G_n$ is a clopen subgroup of $H$ for some $n<\omega$, from which we can see this $G_n$ is also $\alpha$-CLI. Now by (1), we conclude that for all $i\ge n$, $G^i$ is $\alpha$-CLI.

On the other hand, assume all $G^i$ are L-$\alpha$-CLI, and there is an $m<\omega$ such that $G^i$ is $\alpha$-CLI for each $i\ge m$. Then (1) implies that $\prod_{i\ge m}G^i$ is $\alpha$-CLI. Note that $G=G^0\times\dots\times G^{m-1}\times\prod_{i\ge m}G^i$. By Corollary~\ref{GtimesH}, we have
$${\rm rank}(G)=\max\{{\rm rank}(G^0),\dots,{\rm rank}(G^{m-1}),{\rm rank}(\prod_{i\ge m}G^i)\}\le\alpha,$$
i.e., $G$ is L-$\alpha$-CLI.
\end{proof}

In the rest of this article, we will show that the notions of $\alpha$-CLI, L-$\alpha$-CLI, together with ${\rm rank}(G)$, form a proper hierarchy on the class of non-archimedean CLI Polish groups. For this purpose, we will construct groups which are $\alpha$-CLI but not L-$\beta$-CLI for all $\beta<\alpha$, and groups which are L-$\alpha$-CLI but not $\alpha$-CLI, by induction on $\alpha<\omega_1$. We consider the case concerning successor ordinals first.

\begin{corollary}\label{+1}
Let $(G^i)$ be a sequence of non-archimedean CLI Polish groups, $\alpha<\omega_1$, and let $G=\prod_iG^i$. If all $G^i$ are L-$\alpha$-CLI but not $\alpha$-CLI, then $G$ is $(\alpha+1)$-CLI but not L-$\alpha$-CLI.
\end{corollary}

\begin{proof}
Note that all $G^i$ are $(\alpha+1)$-CLI but not $\alpha$-CLI.
\end{proof}

From Theorem~\ref{0-1} and \cite[Theorem 1.1]{GX}, a non-archimedean CLI Polish group $G$ is $1$-CLI iff $G$ is isomorphic to a closed subgroup of a product $\prod_iG^i$, where each $G^i$ is L-$0$-CLI. However, we do not know whether the following generalization of Corollary~\ref*{+1} is true:

\begin{question}
    Let $0<\alpha<\omega_1$, and let $G$ be an $(\alpha+1)$-CLI group. Can we find a sequence of L-$\alpha$-CLI groups $G^i$ such that $G$ is isomorphic to a closed subgroup of $\prod_iG^i$?
\end{question}

Let $G$ and $\Lambda$ be two groups. Recall that the {\it wreath product} $\Lambda\wr G$ is the set $\Lambda\times G^\Lambda$ with the following group operation: given $(\hat\lambda,\hat\chi),(\tilde\lambda,\tilde\chi)\in\Lambda\times G^\Lambda$, we have
$$(\hat\lambda,\hat\chi)(\tilde\lambda,\tilde\chi)=(\hat\lambda\tilde\lambda,\chi),$$
where $\chi(\lambda)=\hat\chi(\lambda)\tilde\chi(\hat\lambda^{-1}\lambda)$ for each $\lambda\in\Lambda$. If $\Lambda$ is countable discrete and $G$ is Polish, then $\Lambda\wr G$ equipped with the product topology of $\Lambda\times G^\Lambda$ is also a Polish group.

\begin{theorem}\label{wr}
Let $G$ be a non-archimedean CLI Polish group, $\Lambda$ an infinite countable discrete group, and $\alpha<\omega_1$. If $G$ is $(\alpha+1)$-CLI but not $\alpha$-CLI, then $\Lambda\wr G$ is L-$(\alpha+1)$-CLI but not $(\alpha+1)$-CLI.
\end{theorem}

\begin{proof}
For brevity, we denote $\Lambda\wr G$ by $H$. Let $\lambda_i,\,i<\omega$ be an enumeration of $\Lambda$ without repetition. Note that the underlying space of $\Lambda\wr G$ is $\Lambda\times G^\Lambda$. For $(\lambda,\chi)\in\Lambda\wr G$, put $\pi_\Lambda(\lambda,\chi)=\lambda$ and $\pi_G^i(\lambda,\chi)=\chi(\lambda_i)$.

Let $\mathcal G=(G_n)\in{\rm dgnb}(G)$. Since $G$ is not $\alpha$-CLI, $G\ne\{1_G\}$. Without loss of generality, we can assume that $G\ne G_1$. Put $H_0=H$ and 
$$H_{n+1}=\{(1_\Lambda,\chi):\chi\in G^\Lambda\wedge\forall i<n\,(\chi(\lambda_i)\in G_n)\}$$
for $n<\omega$. Then $\mathcal H=(H_n)\in{\rm dgnb}(H)$. Note that the open subgroup $H_1=\{1_\Lambda\}\times G^\Lambda$ is topologically isomorphic to $G^\omega$. By Theorem~\ref{prod}, $H_1$ is $(\alpha+1)$-CLI. So we can see that $H$ is L-$(\alpha+1)$-CLI from Theorem~\ref*{locally}.

Since $G$ is $(\alpha+1)$-CLI but not $\alpha$-CLI, $\omega\cdot\alpha<\rho(\mathcal G)\le\omega\cdot(\alpha+1)$. By Lemma~\ref{rank}, there exist $1\le m,k<\omega$ such that $\rho^k(\mathcal G)=\omega\cdot\alpha+m$. To see that $H$ is not $(\alpha+1)$-CLI, we will show that $\rho^{k+1}(\mathcal H)>\omega\cdot(\alpha+1)$ as follows.

For any $(\lambda_l,\hat\chi)\in H$ and $(1_\Lambda,\tilde\chi)\in H_{k+1}$, note that $(\lambda_l,\hat\chi)(1_\Lambda,\tilde\chi)=(\lambda_l,\hat\chi\tilde\chi_l)$, where $\tilde\chi_l(\lambda)=\tilde\chi(\lambda_l^{-1}\lambda)$ for $\lambda\in\Lambda$. It follows that
$$(\lambda_l,\hat\chi)H_{k+1}=\{(\lambda_l,\chi):\chi\in G^\Lambda\wedge\forall i<k\,(\chi(\lambda_l\lambda_i)\in\hat\chi(\lambda_l\lambda_i)G_k)\}.$$
There is a unique $l_i<\omega$ such that $\lambda_{l_i}=\lambda_l\lambda_i$ for $i<k$. It is clear that $\pi_\Lambda((\lambda_l,\hat\chi)H_{k+1})=\{\lambda_l\}$ and for $j<\omega$,
$$\pi_G^j((\lambda_l,\hat\chi)H_{k+1})=\left\{\begin{array}{ll}\hat\chi(\lambda_{l_i})G_k, & j=l_i,i<k,\cr G, & \mbox{otherwise.}\end{array}\right.$$
Let $m_l=\max\{l_i:i<k\}$. Then it is clear that $m_l\ge k-1$ for $l<\omega$ and $\sup\{m_l:l<\omega\}=\omega$. Note that
$\pi_G^{m_l}(H_{m_l+1}\cdot(\lambda_l,\hat\chi)H_{k+1})=G/G_k$ is not a singleton, and so $(m_l+1,H_{m_l+1}\cdot(\lambda_l,\hat\chi)H_{k+1})\in T^{H/H_{k+1}}_\mathcal H$. Also note that
$$(T^{H/H_{k+1}}_\mathcal H)_{(m_l+2,H_{m_l+2}\cdot(\lambda_l,\hat\chi)H_{k+1})}\cong T^{H_{m_l+2}\cdot(\lambda_l,\hat\chi)H_{k+1}}_{\mathcal H'}, \textrm{ and}$$
$$(T^{G/G_k}_\mathcal G)_{(m_l+1,G_{m_l+1}\cdot gG_k)}\cong T^{G_{m_l+1}\cdot gG_k}_{\mathcal G'},$$
where $\mathcal H'=(H_{n+m_l+2})$ and $\mathcal G'=(G_{n+m_l+1})$ for $n<\omega$. Define $\theta:H/H_{k+1}\to G/G_k$ as $\theta(C)=\pi_G^{m_l}(C)$.
Now we restrict $\theta$ as a function that maps $H_{m_l+2}\cdot(\lambda_l,\hat\chi)H_{k+1}$ onto $G_{m_l+1}\cdot\hat\chi(\lambda_{m_l})G_k$. Then apply Lemma~\ref*{surjection-tree} and Proposition~\ref*{L_k(T)} to obtain the following: 
$$\begin{array}{ll}&\rho((T^{H/H_{k+1}}_\mathcal H)_{(m_l+1,H_{m_l+1}\cdot(\lambda_l,\bar\chi)H_{k+1})})\cr
\ge &\sup\{\rho((T^{H/H_{k+1}}_\mathcal H)_{(m_l+2,H_{m_l+2}\cdot(\lambda_l,\hat\chi)H_{k+1})}):\forall i<m_l\,(\bar\chi(\lambda_i)=\hat\chi(\lambda_i))\}+1\cr
\ge &\sup\{\rho((T^{G/G_k}_\mathcal G)_{(m_l+1,G_{m_l+1}\cdot gG_k)}):g\in G\}+1\cr
\ge &\omega(\rho((T^{G/G_k}_\mathcal G))+1\cr
= &\omega\cdot\alpha+1.\end{array}$$
This implies $\rho(T^{H/H_{k+1}}_\mathcal H)\ge\omega\cdot\alpha+m_l+2$ for all $l<\omega$, and hence
$$\rho^{k+1}(\mathcal H)=\rho(T^{H/H_{k+1}}_\mathcal H)\ge\omega\cdot\alpha+\omega=\omega\cdot(\alpha+1).$$
Since $\rho^{k+1}(\mathcal H)$ is a successor ordinal, we have $\rho^{k+1}(\mathcal H)>\omega\cdot(\alpha+1)$.
\end{proof}

Now we turn to the case concerning limit ordinals. To do this, we need two lemmas.

\begin{lemma}
Let $(G^i)$ be a sequence of Polish groups, and let $H^i$ be an open subgroup of $G^i$ for each $i<\omega$. Suppose $H=\prod_iH^i$ and
$$G=\{(g_i)\in\prod_iG^i:\forall^\infty i\,(g_i\in H^i)\}$$
is equipped with the topology $\tau$ generated by the sets of the form $(g_i)U$ for $(g_i)\in G$ and $U$ open in $H$. Then $(G,\tau)$ is a Polish group and $\tau$ is the unique group topology on $G$ making $H$ an open subgroup of $G$.
\end{lemma}

\begin{proof}
For each $(g_i)\in G$, the subspace $(g_i)H$ of $(G,\tau)$ is homeomorphic to $H$, so it is Polish. Let $D^i$ be a subset of $G^i$ which meets every coset of $H^i$ at exactly one point. Note that $D^i$ is countable for all $i<\omega$. We define
$$D=\{(g_i)\in\prod_iD^i:\forall^\infty i\,(g_i=1_{G^i})\}.$$
It is clear that $G/H=\{(g_i)H:(g_i)\in D\}$ is countable. So $(G,\tau)$ is a sum of countably many Polish spaces, thus is a Polish space.

For $(g_i),(h_i)\in G$ and $U$ open in $H$ with $1_H\in H$, there exists an $m<\omega$ such that $g_i,h_i\in H^i$ for $i>m$ and $U^0\times\dots\times U^m\times\prod_{i>m}H^i\subseteq U$, where $U^i$ is an open subset of $H^i$ with $1_{H^i}\in U^i$ for each $i\le m$. We can find open neighborhoods $V^i$ and $W^i$ of $1_{H^i}$ with $(g_iV^i)(h_iW^i)^{-1}\subseteq g_ih_i^{-1}U^i$ for $i\le m$. Now let $V=V^0\times\dots\times V^m\times\prod_{i\ge m}H^i$ and $W=W^0\times\dots\times W^m\times\prod_{i\ge m}H^i$, then $V$ and $W$ are open neighborhoods of $1_H$, and $((g_i)V)((h_i)W)^{-1}\subseteq(g_ih_i^{-1})U$. So $(G,\tau)$ is a Polish group.

Finally, suppose $\tau'$ is another group topology on $G$ such that $H$ is an open subgroup of $G$. Then for each $(g_i)\in G$, the subspace $(g_i)H$ of $(G,\tau')$ is homeomorphic to $H$, so the restrictions of $\tau$ and $\tau'$ on $(g_i)H$ are the same. Hence, $\tau=\tau'$.
\end{proof}

\begin{lemma}\label{limit}
Let $(G^i)$ be a sequence of non-archimedean CLI Polish groups, $\mathcal G^i=(G^i_n)\in{\rm dgnb}(G^i)$ for each $i<\omega$, and let $0<\alpha<\omega_1$. Suppose
$$\sup\{\rho^1(\mathcal G^i):i<\omega\}=\omega\cdot\alpha,\textrm{ and}$$
$$G=\{(g_i)\in\prod_iG^i:\forall^\infty i\,(g_i\in G^i_1)\}$$
is equipped with the unique group topology making $\prod_iG^i_1$ an open subgroup of $G$. Then $G$ is not $\alpha$-CLI.
\end{lemma}

\begin{proof}
Put $G_0=G$ and for $n\ge 1$, let
$$G_n=\prod_{i<n-1}G^i_n\times\prod_{i\ge n-1}G^i_1.$$
It is clear that $G_1=\prod_iG^i_1$ and $\mathcal G=(G_n)\in{\rm dgnb}(G)$.

Given $j<\omega$, define $\theta:G/G_1\to G^j/G^j_1$ as $\theta((g_i)G_1)=g_jG^j_1$ for $(g_i)\in G$. Applying Lemma~\ref{surjection-tree} to the restriction of $\theta$ as in the proof of Theorem~\ref{wr}, we have
$$\begin{array}{ll}\rho(T^{G/G_1}_\mathcal G) &=\sup\{\rho((T^{G/G_1}_\mathcal G)_{(1,G_1\cdot(g_i)G_1)}):(g_i)\in G\}+1\cr
&\ge\sup\{\rho((T^{G^j/G^j_1}_{\mathcal G^j})_{(1,G^j_1\cdot gG^j_1)}):g\in G^j\}+1\cr
&=\rho(T^{G^j/G^j_1}_{\mathcal G^j})=\rho^1(\mathcal G^j).\end{array}$$

Therefore,
$$\rho^1(\mathcal G)=\rho(T^{G/G_1}_\mathcal G)\ge\sup\{\rho^1(\mathcal G^j):j<\omega\}=\omega\cdot\alpha.$$
Since $\rho^1(\mathcal G)$ is a non-limit ordinal and $\alpha>0$, we have $\rho^1(\mathcal G)>\omega\cdot\alpha$. So $G$ is not $\alpha$-CLI.
\end{proof}

Finally, we complete the construction in the following theorem.

\begin{theorem}
For any $\alpha<\omega_1$, there exist two non-archimedean CLI Polish groups $G_\alpha$ and $H_\alpha$ with ${\rm rank}(G_\alpha)={\rm rank}(H_\alpha)=\alpha$ such that $H_\alpha$ is $\alpha$-CLI and $G_\alpha$ is L-$\alpha$-CLI but not $\alpha$-CLI.
\end{theorem}

\begin{proof}
We construct $G_\alpha$ and $H_\alpha$ by induction on $\alpha$. From Corollary~\ref{+1} and Theorem~\ref{wr}, we only need to consider the case where $\alpha$ is a limit ordinal.

Let $(\alpha_i)$ be a sequence of ordinals less than $\alpha$ with $\sup\{\alpha_i:i<\omega\}=\alpha$. By the inductive hypothesis, we can find a non-archimedean CLI Polish group $G^i$ for each $i<\omega$ such that $G^i$ is L-$\alpha_i$-CLI but not $\alpha_i$-CLI. It is clear that ${\rm rank}(G^i)=\alpha_i<\alpha$.

Put $H_\alpha=\prod_iG^i$. Theorem~\ref{prod}.(1) implies that $H_\alpha$ is $\alpha$-CLI. By Theorem~\ref{subgroup}, ${\rm rank}(H_\alpha)\ge{\rm rank}(G^i)$ for each $i<\omega$. So ${\rm rank}(H_\alpha)=\alpha$.

For $i<\omega$, let $\mathcal G^i=(G^i_n)\in{\rm dgnb}(G^i)$. By Lemma~\ref{rank}, there exist $0<k_i<\omega$ such that $\omega(\rho^{k_i}(\mathcal G^i))=\omega\cdot\alpha_i$. Put $H^i_0=G^i$ and $H^i_{n+1}=G^i_{n+k_i}$ for $n<\omega$. Then $\mathcal H^i=(H^i_n)\in{\rm dgnb}(G^i)$ and $H^i_1=G^i_{k_i}$. By Lemma~\ref{two dgnb}, we have
$$\omega(\rho^1(\mathcal H^i))=\omega(\rho(T^{G^i/H^i_1}_{\mathcal H^i}))=\omega(\rho(T^{G^i/G^i_{k_i}}_{\mathcal G^i}))=\omega(\rho^{k_i}(\mathcal G^i))=\omega\cdot\alpha_i.$$
So $\sup\{\rho^1(\mathcal H^i):i<\omega\}=\omega\cdot\alpha$. Now we let
$$G_\alpha=\{(g_i)\in\prod_iG^i:\forall^\infty i\,(g_i\in G^i_{k_i})\}$$
be equipped with the unique group topology making $\prod_iG^i_{k_i}$ an open subgroup of $G_\alpha$. By Lemma~\ref{limit}, $G_\alpha$ is not $\alpha$-CLI. It is clear that the open subgroup $\prod_iG^i_{k_i}$ is $\alpha$-CLI, so $G_\alpha$ is L-$\alpha$-CLI, and hence ${\rm rank}(G_\alpha)=\alpha$.
\end{proof}

\subsection*{Acknowledgements}
We are grateful to Su Gao and Victor Hugo Ya\~{n}ez for carefully reading the manuscript, and giving much useful advice on grammar and language. We also thank the anonymous referee for helpful suggestions in both the context of definitions and language.

\end{document}